\documentclass[12pt, oneside, a4paper]{article}
\usepackage{graphicx}
\usepackage{tikz}
\usepackage{cite}
\usepackage{hyperref}
\usepackage{amsfonts}
\usepackage{fancyhdr,amsthm,amsmath,amssymb,mathrsfs}
\usetikzlibrary{intersections}

\newtheorem{thm}{\bf Theorem}[section]
\newtheorem{lem}[thm]{Lemma}

\newtheorem{cor}[thm]{Corollary}
\theoremstyle{definition}

\newtheorem{rem}[thm]{Remark}
\newtheorem{eg}[thm]{Example}
\numberwithin{equation}{section}
\newenvironment{prf}{\noindent {\bf Proof.}\rm}{\qed}
\newenvironment{prff}{\noindent {\bf Proof of Theorem \ref{main thm}.}\rm}{\qed}
\newcommand{\pth}{{\rm Path}}
\newcommand{\In}{{\rm Inv}(\Gamma)}
\newcommand{\Ini}{{\rm Inv}(\Gamma_i)}
\newcommand{\Inx}{{\rm Inv}(\Gamma_1)}
\newcommand{\Iny}{{\rm Inv}(\Gamma_2)}
\newcommand{\al}{\alpha}

\newcommand{\G}{\Gamma}

\newcommand{\DOI}[1]{DOI: \href{https://doi.org/#1}{#1}}
\newcommand{\arxiv}[1]{arXiv: \href{https://doi.org/#1}{#1}}

\begin{document}
\title{\Large\bf Distributivity in congruence lattices of graph inverse semigroups\thanks{Supported by the National Natural Science
Foundation of China (Grant No. 11771212) and the Postgraduate Research $\&$ Practice Innovation Program of Jiangsu Province.}}
\author{Yongle Luo$^{a, b}$, Zhengpan Wang$^{b}$ and Jiaqun Wei$^{a}$\thanks{Corresponding author} \\
{\small\em $^{a}$School of Mathematics Sciences, Nanjing Normal University, $210023$ Nanjing, China} \\
{\small\em $^{b}$School of Mathematics and Statistics, Southwest University, $400715$ Chongqing, China} \\
{\small\em E-mail: yongleluo@nnu.edu.cn, \; zpwang@swu.edu.cn, \; weijiaqun@njnu.edu.cn} }
\date{}
\maketitle

%%%%%%%%%%%%%%%%%%%%%%%%%%%%%%%%%%%%%%%%%%%%%%%%%%%%%%%%%%%%%%%%%%%%%%%%%%%%%%%%%%%%%%%
%--------------------------------------------------------------------------------------

\begin{abstract}
Let $\G$ be a directed graph and $\In$ be the graph inverse semigroup of $\G$. Luo and Wang \cite{Luo} showed that the congruence lattice $\mathscr{C}(\In)$ of any graph inverse semigroup $\In$ is upper semimodular, but not lower semimodular in general. Anagnostopoulou-Merkouri, Mesyan and Mitchell characterized the directed graph $\G$ for which $\mathscr{C}(\In)$ is lower semimodular \cite{Anag}. In the present paper, we show that the lower semimodularity, modularity and distributivity in the congruence lattice $\mathscr{C}(\In)$ of any graph inverse semigroup $\In$ are equivalent.\medskip


\end{abstract}

\section{Introduction and main result}

The earliest research on graph inverse semigroups can be traced back to 1975. Ash and Hall in \cite{AshHall} introduced these semigroups to study the nonzero $\mathscr{J}$-classes of inverse semigroups. Graph inverse semigroups are closely related to the study of Leavitt path algebras \cite{AbramsPino,Meakin} and Cuntz inverse semigroups in the $C^\ast$-algebras \cite{KumjianPaskRaeburn,Paterson}.

A {\it directed graph} $\G=(V, E, s, r)$ is a quadruple consisting of two sets $V, E$ and two functions $s, r : E \rightarrow V$. The elements of $V$ are called {\it vertices} and the elements of $E$ are called {\it edges}. For each edge $e$ in $E$, $s(e)$ is the {\it source} of $e$ and $r(e)$ is the {\it range} of $e$. In this case, we say that the {\it direction} of $e$ is from $s(e)$ to $r(e)$. For any set $X$, we denote by $|X|$ the {\it cardinality} of $X$. For any vertex $v$, $|s^{-1}(v)|$ is called the {\it index} of $v$. Where there is no possibility for confusion, we may write $\G=(V, E, s, r)$ simply as $\G$.  A directed graph $\G$ is {\it finite} if both $V$ and $E$ are finite.

A {\it nontrivial path} in $\G$ is a sequence $\al := e_1 \cdots e_n$ of edges such that $r(e_i) = s(e_{i+1})$ for $i = 1, \cdots, n-1$. In this case, $s(\al):= s(e_1)$ is the {\it source} of $\al$, $r(\al):= r(e_n)$ is the {\it range} of $\al$. We denote by $u(\alpha)$ the set $\{s(e_i): i=1, \cdots, n\}$. A nontrivial path $\al$ is called a {\it cycle} if $s(\al) = r(\al)$ and $s(e_i) \neq s(e_j)$ whenever $i \neq j$. Given a subset $V_1$ of $V$, we denote by $C({V_1})$ the set of all cycles $c$ with $u(c) \subseteq V_1$. We say that $\G$ is {\it acyclic} if it contains no cycles, that is, $C(V)=\emptyset$. Moreover, we also view a vertex $v$ of $V$ as a {\it trivial path} with $s(v) = r(v) = v$. In this case, we naturally take $u(v) = \emptyset$, since it has no edges. Denote by $\pth(\G)$ the set of all paths in $\G$.

The {\it graph inverse semigroup} $\In$ of $\G$ is the semigroup with zero generated by $V$, $E$ and $E^\ast$, where $E^\ast$ is the set $\{e^\ast : e \in E\}$, satisfying the following relations for all $u, v \in V$ and $e, f \in E$:
\begin{description}
\item[{\rm(GI1)}] $u v = \delta_{u,v} u$;
\item[{\rm(GI2)}] $s(e) e = e r(e) = e$;
\item[{\rm(GI3)}] $r(e) e^\ast = e^\ast s(e) = e^\ast$;
\item[{\rm(CK1)}] $f^\ast e = \delta_{f,e} r(e)$.
\end{description}
The symbol $\delta$ is the Kronecker delta.

Certainly, $\In$ is finite if and only if $\G$ is finite and acyclic. One can see that every nonzero element in $\In$ can be uniquely written as $\alpha\beta^{*}$ for some $\alpha, \beta$ in $\pth(\G)$ with $r(\alpha) = r(\beta)$. The semigroup $\In$ is an inverse semigroup with semilattice $\mathcal{E}(\In)$ of idempotents, where
$$ \mathcal{E}(\In)=\{\alpha\alpha^{*}:\alpha \in \pth(\G)\}\cup \{0\}.$$
Moreover, the multiplication on $\In$ is also defined as follows:
\begin{equation} \label{multiplication of Inv}
(\alpha\beta^*)(\zeta\eta^*) = {\begin{cases}
(\alpha\xi)\eta^* & \mbox{if } \zeta=\beta\xi ~ \mbox{for}~\mbox{some}~\xi\in \pth(\G), \\
\alpha(\eta\xi)^*  & \mbox{if } \beta=\zeta\xi~\mbox{for}~\mbox{some}~\xi\in \pth(\G), \\
0 & \mbox{otherwise}.
\end{cases}}
\end{equation}

There are many papers about congruences and congruence lattices of graph inverse semigroups. Mesyan and Mitchell \cite{MesyanMitchell} showed that the Green's $\mathscr{H}$-relation on any graph inverse semigroup $\In$ is trivial. Thus every graph inverse semigroup is combinatorial. For a finite combinatorial inverse semigroup, Jones \cite{Jones} proved that the congruence lattice is $M$-symmetric, hence upper semimodular. Wang \cite{Wang} described the congruences of $\In$ in terms of certain sets of vertices and integer-valued functions on the cycles in $\G$. Based on this description, Wang showed that each graph inverse semigroup $\In$ of a finite directed graph $\G$ is congruence Noetherian. Luo and Wang \cite{Luo} proved that the congruence lattice $\mathscr{C}(\In)$ of any graph inverse semigroup $\In$ is upper semimodular, but not lower semimodular in general. In \cite{Anag}, Anagnostopoulou-Merkouri, Mesyan and Mitchell provided a characterization of the directed graph $\G$ for which $\mathscr{C}(\In)$ is lower semimodular.

Based on the results of \cite{Anag}, \cite{Luo} and \cite{Wang}, we further show that the lower semimodularity, modularity and distributivity in the congruence lattice $\mathscr{C}(\In)$ of any graph inverse semigroup $\In$ are equivalent, which extends the result in \cite[Corollary 2.4]{Anag}. The main result in the present paper is as follows.

\begin{thm}\label{main thm}
Let $\G$ be a directed graph and $\mathscr{C}(\In)$ be the congruence lattice of the graph inverse semigroup $\In$. Then the following statements are equivalent.
\begin{description}
\item[(\romannumeral1)] The lattice $\mathscr{C}(\In)$ is lower semimodular;
\item[(\romannumeral2)] The lattice $\mathscr{C}(\In)$ is modular;
\item[(\romannumeral3)] The lattice $\mathscr{C}(\In)$ is distributive.
\end{description}
\end{thm}

\section{Preliminaries}

\subsection{\textit{Lattices}}

In this subsection, we recall some fundamental materials on lattices that can be found in the standard lattice theory book \cite{Gratzer}.

Let $(L, \leq)$ be a partially ordered set. For any $a, b \in L$, we say that $a$ {\it covers} $b$, and write $a\succ b$, if $a > b$ and there is no $x \in L$ such that $a > x > b$. A partially ordered set $(L, \leq)$ is called a {\it lattice} if, for all $a, b \in L$, there exists a least upper bound $a\vee b\in L$, called the {\it join} of $a$ and $b$, and a greatest lower bound $a\wedge b \in L$, called the {\it meet} of $a$ and $b$.

Let $L=(L,\leq, \vee, \wedge)$ be a lattice. Recall that $L$ is {\it distributive} if $(a\vee b)\wedge c =(a\wedge c)\vee (b\wedge c)$, for all $a, b, c \in L$. Moreover, $L$ is {\it modular} if $a \leq c$ implies that $a\vee (b\wedge c)= (a\vee b)\wedge c$, for all $a, b, c \in L$. The lattice $L$ is {\it upper semimodular} if $a, b\succ a\wedge b$ implies that $a \vee b \succ a, b$, for all $a, b \in L$. Dually, $L$ is {\it lower semimodular} if
$a\vee b \succ a, b$ implies that $a, b\succ a\wedge b$, for all $a, b \in L$.

It is clear that every distributive lattice is modular. Every modular lattice is both upper semimodular and lower semimodular, but the converse does not hold for some infinite lattices.

\begin{figure}[htbp]
\centering
\begin{tikzpicture}[scale=0.7]
\fill(0,2.5)circle(1pt);
\fill (1,1)circle(1pt);
\fill (2,2)circle(1pt);
\fill (2,3)circle(1pt);
\fill (1,4)circle(1pt);
\draw (1,1)--(2,2);
\draw (2,2)--(2,3);
\draw (1,4)--(2,3);
\draw (0,2.5)--(1,4);
\draw (0,2.5)--(1,1);
\node[below] at (1,1) {$\mathfrak{N}_5$};

\fill(5,1)circle(1pt);
\fill (5,2.5)circle(1pt);
\fill (5,4)circle(1pt);
\fill (4,2.5)circle(1pt);
\fill (6,2.5)circle(1pt);
\draw (5,1)--(5,2.5);
\draw (5,2.5)--(5,4);
\draw (4,2.5)--(5,4);
\draw (6,2.5)--(5,4);
\draw (5,1)--(4,2.5);
\draw (5,1)--(6,2.5);
\node[below] at (5,1) {$\mathfrak{M}_3$};
\end{tikzpicture}
\caption{\footnotesize The lattices $\mathfrak{N}_5$ and $\mathfrak{M}_3$.}\label{Figure 1}
\end{figure}
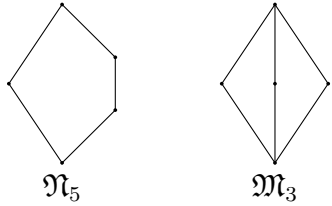

A sublattice $L_1$ of a lattice $L$ is called a {\it pentagon}, respectively a {\it diamond}, if $L_1$ is isomorphic to $\mathfrak{N}_5$, respectively to $\mathfrak{M}_3$, shown in Figure \ref{Figure 1}. The following two lemmas give some characterizations of a lattice to be distributive or modular.

\begin{lem}\cite [Theorem 101]{Gratzer})\label{diagrams of distri}
A lattice $L$ is distributive if and only if it does not contain a pentagon or a diamond.
\end{lem}

\begin{lem}\cite[Theorem 102]{Gratzer})\label{diagrams of modular}
A lattice $L$ is modular if and only if it does not contain a pentagon.
\end{lem}

\subsection{\textit{Graphs}}

Let $\G=(V, E, s, r)$ be a directed graph. For any path $\alpha \in \pth(\G)$, we say that the path $\alpha$ is {\it from $v_1$ to $v_2$} if $s(\alpha) = v_1$ and $r(\alpha)= v_2$. A directed graph $\G$ is {\it strongly connected} if, for all vertices $v_{1}, v_2 \in V$, there exists both a path from $v_1$ to $v_2$ and a path from $v_2$ to $v_1$. Moreover, $\G$ is {\it unilaterally connected} if, for all vertices $v_{1}, v_2 \in V$, there exists either a path from $v_1$ to $v_2$ or a path from $v_2$ to $v_1$.

For any edge $e\in E$, when we do not consider its direction, we say that the unordered pair $(s(e), r(e))=(r(e), s(e))$ is the {\it undirected edge} corresponding to $e$. Denote by $\overline{E}$ the set $\{(s(e), r(e)): e\in E\}$. Then the pair $\overline{\G}=(V, \overline{E})$ is called the {\it underlying graph} of $\G$. Roughly speaking, $\overline{\G}$ is the undirected graph obtained from $\G$ by forgetting the direction of each edge. A {\it nontrivial undirected path} in $\overline{\G}$ is a sequence $\overline{\al}:= (v_1, v_2)(v_2, v_3)\cdots(v_{n-1}, v_n)$ of undirected edges. In this case, we say that $v_1$ and $v_n$ are connected by $\overline{\al}$.

A directed graph $\G$ is {\it weakly connected} if its underlying graph $\overline{\G}$ is connected (i.e., any two distinct vertices are connected by some nontrivial undirected path). A {\it subgraph} of $\G$ is a directed graph $\G_i = (V_{i}, E_{i}, s_{i}, r_{i})$ such that $V_i \subseteq V$, $E_i \subseteq E$ and the restrictions $s\mid_{E_{i}}$, $r\mid_{E_{i}}$ of $s, r$ to $E_{i}$ are respectively equal to $s_{i}, r_{i}$.  A subgraph $\G_{i}$ of $\G$ is called a {\it weakly connected component} if $\G_{i}$ is the maximal weakly connected subgraph of $\G$.

In addition, we need the following notions, which are used to describe the lower semimodularity of the congruence lattice $\mathscr{C}(\In)$ in \cite{Anag}.

Given a directed graph $\G$, there is a natural {\it preorder} $\trianglerighteq$ on $V$: for any $v_1, v_2 \in V$, $v_1 \trianglerighteq v_2$ if there exists a path from $v_1$ to $v_2$. We say that $v \in V$ is a {\it forked vertex} if there exist distinct edges $e, f \in s^{-1}(v)$ such that the following properties hold:
\begin{description}
\item[(\romannumeral1)] $r(g)\ntrianglerighteq r(e)$ for all $g\in s^{-1}(v) \setminus \{e\}$;
\item[(\romannumeral2)] $r(g)\ntrianglerighteq r(f)$ for all $g\in s^{-1}(v) \setminus \{g\}$.
\end{description}

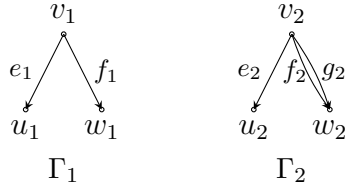
\begin{figure}[htbp]
\centering
\begin{tikzpicture}
\draw (1.5,1) circle(0.03);
\draw (1,0) circle(0.03);
\draw (2,0) circle(0.03);
\node[below] at (1,0) {$u_1$}; \node[below] at (2,0) {$w_1$}; \node[above] at (1.5,1) {$v_1$}; \node[below] at (1.5,-0.5) {\small {$\G_1$}};
\draw[-stealth] (1.5,1) to (1,0); \draw[-stealth] (1.5,1) to (2,0);
\node[left] at (1.25,0.5) {\footnotesize {$e_1$}}; \node[right] at (1.75,0.5) {\footnotesize {$f_1$}};

\draw (4.5,1) circle(0.03); \draw (4,0) circle(0.03); \draw (5,0) circle(0.03);
\node[below] at (4,0) {$u_2$}; \node[below] at (5,0) {$w_2$}; \node[above] at (4.5,1) {$v_2$}; \node[below] at (4.5,-0.5) {\small {$\G_2$}};

\node[left] at (4.25,0.5) {\footnotesize {$e_2$}}; \node[right] at (4.75,0.5) {\footnotesize {$g_2$}}; \node[left] at (4.85,0.5) {\footnotesize {$f_2$}};
\draw[-stealth] (4.5,1) to (4,0); \draw[-stealth] (4.5,1) to [bend right = 8] (5,0); \draw[-stealth] (4.5,1) to [bend left = 8] (5,0);
\end{tikzpicture}	
\caption{\footnotesize There are two directed graphs. The directed graph $\G_1$ has a forked vertex $v_1$. But the directed graph $\G_2$ has no forked vertices.}\label{Figure 2}
\end{figure}

For instance, there are two directed graphs shown in Figure \ref{Figure 2}. The directed graph $\G_1$ has a forked vertex $v_1$. But $\G_2$ has no forked vertices.

\subsection{\textit{Congruences and congruence triples}}

We need Wang's results in \cite{Wang}, which proved that there is a one-to-one correspondence between congruences on a graph inverse semigroup and so-called congruence triples. For convenience, we briefly recall the relevant concepts as follows (see \cite{Luo}).

Given a directed graph $\G$, a subset $H$ of $V$ is {\it hereditary} if $s(e) \in H$ always implies $r(e) \in H$ for all $e \in E$. Let $H$ be a hereditary subset of $V$. For any $v \in V$, we call $|\{e \in E: s(e) = v, r(e) \notin H \}|$ the {\it index of $v$ relative to $H$}. Let $W$ be a subset of vertices with index 1 relative to $H$. Certainly, $H$ and $W$ are two disjoint subsets. We call $f: C(V)\rightarrow \mathbb{Z}^+ \cup \{\infty\}$ a {\it cycle function} if $f$ satisfies the following conditions:
\begin{description}
\item[(\romannumeral1)] $f(c)=1$, if $c \in C(H)$;
\item[(\romannumeral2)] $f(c) = f(c')$, if $c$ and $c'$ are cyclic permutations of each other;
\item[(\romannumeral3)] $f(c)=\infty$, if $c\notin C(H)$ and $c\notin C(W)$.
\end{description}

A {\it congruence triple of $\G$}, written by $(H, W, f)$, consists of a hereditary subset $H$, a subset $W$ of the set of vertices with index 1 relative to $H$, and a cycle function $f$. If $\G$ is acyclic, that is, $C(V) = \emptyset$, then we denote by $(H, W, \emptyset)$ the corresponding congruence triple.

For any $m_1, m_2 \in \mathbb{Z}^+$, we denote by $gcd(m_1, m_2)$ their {\it greatest common divisor} and $lcm(m_1, m_2)$ their {\it least common multiple}. We naturally assume that any element in $\mathbb{Z}^+ \cup \{\infty\}$ divides $\infty$, which means that $gcd(m, \infty)= m$ and $lcm(m, \infty)= \infty$, for all $m\in \mathbb{Z}^+ \cup \{\infty\}$. There is a partial order on $\mathcal{CT}(\G)$: for any $(H_1, W_1, f_1), (H_2, W_2, f_2) \in \mathcal{CT}(\G)$, $(H_1, W_1, f_1) \leq (H_2, W_2, f_2)$ if $H_1 \subseteq H_2$, $W_1 \setminus H_2 \subseteq W_2$ and $f_2(c) \mid f_1(c)$ for any $c \in C(V)$.

Denote by $\mathcal{CT}(\G)$ the set of all congruence triples of $\G$, and denote by $\mathscr{C}(\In)$ the set of all congruences on $\In$. The following lemma shows that $\mathcal{CT}(\G)$ is a lattice  isomorphic to $\mathscr{C}(\In)$.

\begin{lem}\cite[Proposition 1.2]{Luo} \label{isom}
The two partially ordered sets $(\mathscr{C}(\In), \subseteq)$ and $(\mathcal{CT}(\G), \leq )$ are order isomorphic.
\end{lem}

Let $T_1=(H_1, W_1, f_1)$ and $T_2=(H_2, W_2, f_2)$ be any congruence triples in $\mathcal{CT}(\G)$. It is clear that both $H_1\cap H_2$ and $H_1 \cup H_2$ are hereditary subsets. Given a vertex $v\in W_1\cup W_2$, the index of $v$ relative to $H_1 \cup H_2$ is at most 1. We denote by $V_0$ the set of vertices in $(W_1 \cup W_2) \setminus (H_1 \cup H_2)$ with index 0 relative to $H_1 \cup H_2$.

Given a vertex $v \in W_1 \cap W_2$, the index of $v$ relative to both $H_1$ and $H_2$ is 1. If the index of $v$ relative to $H_1 \cup H_2$ is 0, then $v \in V_0$, and the index of $v$ relative to $H_1 \cap H_2$ is 2. If the index of $v$ relative to $H_1 \cup H_2$ is 1, then the index of $v$ relative to $H_1 \cap H_2$ is also 1.

We denote by $X$ the set $(W_1 \cap W_2) \setminus V_0$ and $J$ the set
$$\{v\in{(W_1\cup W_2) \setminus (H_1\cup H_2)}: (\exists \alpha \in \pth(\G)) \;v=s(\alpha), u(\alpha) \subseteq W_1 \cup W_2, r(\alpha)\in V_0\}.$$
Certainly, $V_0 \subseteq J$ since we consider a vertex as a trivial path.

Moreover, we also need the following lemma.
\begin{lem}\cite[Lemmas 2.7 and 2.8]{Luo}\label{meet and join}
Let $T_1=(H_1, W_1, f_1)$ and $T_2=(H_2, W_2, f_2)$ be arbitrary congruence triples in $\mathcal{CT}(\G)$. Then the following statements hold.
\begin{description}
\item[(\romannumeral1)] The meet $T_l = (H_l, W_l, f_l)$ of $T_1$ and $T_2$ can be determined as follows: $H_l = H_1 \cap H_2$, $W_l =(W_1 \cap H_2)\cup (W_2\cap H_1) \cup X$ and $f_l(c) = lcm(f_1(c),f_2(c))$ for $c \in C(V)$.
\item[(\romannumeral2)] The join $T_u = (H_u, W_u, f_u)$ of $T_1$ and $T_2$ can be determined as follows: $H_u = H_1\cup H_2 \cup J$, $W_u = (W_1 \cup W_2) \backslash{H_u}$ and $f_u(c) = gcd(f_1(c), f_2(c))$ for $c \in C(V)$.
\end{description}
\end{lem}

\begin{rem}
Since $H_1 \cap W_1 = H_2 \cap W_2=\emptyset$, the subsets $W_1 \cap H_2$, $W_2 \cap H_1$ and $W_1 \cap W_2$ are disjoint. Moreover, the join $T_u$ satisfies the property that $H_u \cup W_u = H_1 \cup H_2 \cup W_1 \cup W_2$.
\end{rem}

\section{Proof of the main result}
For any disjoint sets $A$ and $B$, we sometimes write their union as $A \uplus B$. The following lemma will be used to prove Theorem \ref{main thm}.

\begin{lem}\label{property of H2 and H3}
Let $\G$ be a directed graph, and let $T_1=(H_1, W_1, f_1)$, $T_2=(H_2, W_2, f_2)$ and $T_3=(H_3, W_3, f_3)$ be congruence triples in $\mathcal{CT}(\G)$ such that $T_1 \wedge T_2 = T_1 \wedge T_3$ and $T_1 \vee T_2 = T_1 \vee T_3$. Then we have $H_2 = H_3$ and $f_2 = f_3$.
\end{lem}

\begin{prf} Denote by $T_l = (H_l, W_l, f_l)$ the meet of $T_1$ and $T_2$, and $T_u = (H_u, W_u, f_u)$ the join of $T_1$ and $T_2$. We know from Lemma $\ref{meet and join}$ that
$$gcd(f_1(c), f_2(c))= gcd(f_1(c), f_3(c))\quad \mbox{and} \quad lcm(f_1(c), f_2(c))= lcm(f_1(c), f_3(c))$$
for any cycle $c \in C(V)$. By the distributivity of $(\mathbb{Z}^{+} \cup \{\infty\}, \mid)$, we have $f_2(c)=f_3(c)$ for any cycle $c \in C(V)$.

Suppose there is a vertex $v \in H_3 \setminus H_2$. It is clear that $v\notin H_1$. Otherwise, it follows from $v\in H_1 \cap H_3 = H_1 \cap H_2$ that $v\in H_2$. Moreover, since $T_1 \vee T_2 =T_1 \vee T_3$, then we get from Lemma \ref{meet and join} again that
$$H_1 \cup W_1 \cup H_2 \cup W_2 = H_1 \cup W_1 \cup H_3 \cup W_3.$$
Now we have
$$v\in H_3 \cap(W_1 \cup W_2)= (H_3 \cap W_1) \cup (H_3 \cap W_2).$$
If $v\in H_3 \cap W_2$, then certainly $v\in W_2 \setminus H_1$. If $v \in H_3\cap W_1$, then we also get $v\in W_2\setminus H_1$, since
$$H_3 \cap W_1 \subseteq W_l \quad \mbox{and} \quad W_l \subseteq  (W_1 \cap H_2) \uplus (W_2 \cap H_1) \uplus (W_1 \cap W_2).$$

We claim that the index of $v$ relative to $H_1 \cup H_2$ is 1. Indeed, it follows from $v\in W_2$ that there exists a unique edge $e_v \in s^{-1}(v)$ such that $r(e_v)\notin H_2$ and the ranges of the rest of edges in $s^{-1}(v)$ are all in $H_2 \subseteq H_1 \cup H_2$. Since $H_3$ is hereditary, we have $r(e_v)\in H_3\setminus H_2$. Similar to the discussion of $v$, it follows $r(e_v)\notin H_1$, that is, $r(e_v) \notin H_1 \cup H_2$. Thus the claim holds.

On the other hand, since $v\in H_3\subseteq H_u$, by the definition of the join of $T_1$ and $T_2$, there must exist a path $\alpha$ starting at $v$
such that $u(\alpha) \subseteq W_1\cup W_2$, $r(\alpha)\in(W_1\cup W_2)\setminus (H_1 \cup H_2)$, and the index of $r(\alpha)$ relative to $H_1 \cup H_2$ is 0. Note that $r(\alpha) \in H_3 \setminus H_2$. Similar to the discussion of $v$, the index of $r(\alpha)$ relative to $H_1 \cup H_2$ is 1, which is a contradiction. Then we proved that $H_3 \setminus H_2 = \emptyset$, that is, $H_3 \subseteq H_2$.

Dually, we also have $H_2 \subseteq H_3$.
\end{prf}\medskip

\begin{lem}\label{distri=modular}
Let $\G$ be a directed graph and $\mathscr{C}(\In)$ be the congruence lattice of the graph inverse semigroup $\In$. Then $\mathscr{C}(\In)$ is modular if and only if it is distributive.
\end{lem}

\begin{prf} By Lemmas \ref{diagrams of distri}, \ref{diagrams of modular} and \ref{isom}, we only need to prove that $\mathcal{CT}(\G)$ does not have a diamond $\mathfrak{M}_3$. Suppose $\mathcal{CT}(\G)$ has a diamond. Without loss of generality, we may assume that there are three different congruence triples $T_i =(H_i, W_i, f_i)$, $i=1, 2, 3$, such that
$$T_1 \wedge T_2 = T_1\wedge T_3 = T_2 \wedge T_3 \quad \mbox{and} \quad T_1\vee T_2 =T_1 \vee T_3 = T_2 \vee T_3.$$
By Lemma \ref{property of H2 and H3}, we have $H_1 = H_2 = H_3$ and $f_1 = f_2 =f_3$. Then we get from Lemma \ref{meet and join} that
$$W_1 \cap W_2 = W_1 \cap W_3 \quad \mbox{and} \quad W_1 \cup W_2 = W_1 \cup W_3.$$
Hence we have $W_2 = W_3$, which contradicts the fact that $T_2 \neq T_3$.
\end{prf}\medskip

Anagnostopoulou-Merkouri, Mesyan and Mitchell \cite{Anag} characterized the directed graph $\G$ for which the congruence lattice $\mathscr{C}(\In)$ is lower semimodular.

\begin{lem}\cite[Theorem 2.3]{Anag}\label{semimodular by AMM}
Let $\G$ be a directed graph. Then $\G$ has no forked vertices if and only if the congruence lattice $\mathscr{C}(\In)$ of $\In$ is lower semimodular.
\end{lem}

Inspired by this result, we have the following lemma.
\begin{lem}\label{no fork then modular}
Let $\G$ be a directed graph. If $\G$ has no forked vertices, then the congruence lattice $\mathscr{C}(\In)$ of $\In$ is modular.
\end{lem}

\begin{prf} Let $\G$ be a directed graph with no forked vertices. Suppose that $\mathscr{C}(\In)$ is not modular. Equivalently, $\mathcal{CT}(\G)$ is not modular. Then it follows from Lemma \ref{diagrams of modular} that $\mathcal{CT}(\G)$ has a pentagon $\mathfrak{N}_5$. Without loss of generality, we may assume that there are congruence triples $T_i =(H_i, W_i, f_i)$, $i=1, 2, 3$, such that
$$ T_1 \wedge T_2 =T_1 \wedge T_3, \quad T_1 \vee T_2 =T_1 \vee T_3 \quad \mbox{and} \quad T_2< T_3.$$
By Lemma \ref{property of H2 and H3}, we have $H_2 = H_3$, $f_2 =f_3$, and $W_2 \subsetneq W_3$.

Denote by $T_l = (H_l, W_l, f_l)$ the meet of $T_1$ and $T_2$, and $T_u = (H_u, W_u, f_u)$ the join of $T_1$ and $T_2$. Let $X_2$ (respectively, $X_3$) be the subset of $W_1 \cap W_2$ (respectively, $W_1 \cap W_3$) such that the index of each vertex relative to $H_l$ is 1. It follows from Lemma \ref{meet and join} that
$$(W_1  \cap H_2) \uplus (W_2 \cap H_1)\uplus X_2 = (W_1 \cap H_3) \uplus (W_3 \cap H_1)\uplus X_3.$$
Since $H_2 =H_3$, we have $(W_2 \cap H_1) \uplus X_2= (W_3 \cap H_1)\uplus X_3$.
For any vertex $v \in X_2$, it follows from
$$X_2 \cap H_1 \subseteq (W_1 \cap W_2) \cap H_1 \subseteq W_1 \cap H_1 = \emptyset$$
that $v\notin H_1$. Thus we have $v\in X_3$, that is, $X_2 \subseteq X_3$. Similarly, we also have $X_3 \subseteq X_2$. So we obtain from $X_2 =X_3$ that $W_2 \cap H_1=W_3 \cap H_1$.

On the other hand, it follows from $T_1\vee T_2 = T_1 \vee T_3$ that
$$H_1 \cup W_1 \cup H_2 \cup W_2 = H_1 \cup W_1 \cup H_3 \cup W_3.$$
Note that $H_2 =H_3$. Denote by
$$A= H_1 \cup W_1 \cup H_2 = H_1 \cup W_1 \cup H_3.$$
Then we have $W_2\setminus A= W_3\setminus A$. It follows from
$$(W_1 \cap W_2)\uplus (W_2 \cap H_1) \uplus (W_2\setminus A)= W_2 \subsetneq W_3 = (W_1 \cap W_3)\uplus (W_3 \cap H_1) \uplus (W_3\setminus A)$$
that $W_1 \cap W_2 \subsetneq W_1 \cap W_3$. Let $B_2$ (respectively, $B_3$) be the subset of $W_1 \cap W_2$ (respectively, $W_1 \cap W_3$) such that the index of each vertex relative to $H_1 \cup H_2$ is 0. We know from
$$X_2 = X_3 \quad \mbox{and} \quad W_1 \cap W_2 = X_2 \uplus B_2 \subsetneq X_3 \uplus B_3 = W_1 \cap W_3$$
that $B_2 \subsetneq B_3$. Certainly, $B_3 \neq \emptyset$, which means that there is a vertex $v_0$ in $B_3$ such that the index of $v_0$ relative to $H_1 \cup H_2$ is 0. In this case, the index of $v_0$ relative to both $H_1$ and $H_3$ is 1, and its index relative to $H_1 \cap H_3$ is 2. Thus there are at least two distinct edges $e, f \in s^{-1}(v_0)$ such that $r(e) \in H_1\setminus H_3$ and $r(f)\in H_3\setminus H_1$. And if there is any other edge $g \in s^{-1}(v_0)$, then $r(g) \in H_1 \cap H_3$. Note that $H_1$, $H_3$ and $H_1 \cap H_3$ are all hereditary subsets. Thus $v_0$ is a forked vertex, which is a contradiction.
\end{prf}\medskip

\begin{prff}
It is clear that $(\mathbf{\romannumeral3})\Rightarrow (\mathbf{\romannumeral2}) \Rightarrow (\mathbf{\romannumeral1})$. The proof of $(\mathbf{\romannumeral1}) \Rightarrow (\mathbf{\romannumeral3})$ follows from Lemmas \ref{distri=modular}, \ref{semimodular by AMM} and \ref{no fork then modular}.
\end{prff}\medskip

As corollaries, we have the following interesting results.

\begin{cor}
Let $\G$ be a directed graph with $|s^{-1} (v)| \leq 1$ for any vertex $v$ in $V$. Then the congruence lattice $\mathscr{C}(\In)$ of $\In$ is distributive.
\end{cor}

\begin{prf} Since the index of every forked vertex is at least 2, the directed graph $\G$ has no forked vertices. Thus the lattice $\mathscr{C}(\In)$ is distributive by Theorem \ref{main thm} and Lemma \ref{semimodular by AMM}.
\end{prf}\medskip

\begin{cor}
Let $\Gamma$ be a directed graph. If $\G$ is unilaterally connected, in particular, strongly connected, then the congruence lattice $\mathscr{C}(\In)$ of $\In$ is distributive.
\end{cor}

\begin{prf} By Theorem \ref{main thm} and Lemma \ref{semimodular by AMM}, it is enough to show that $\G$ has no forked vertices. Suppose that $v$ is a forked vertex in $V$. For any distinct edges $e, f\in s^{-1}(v)$, it follows from the unilateral connectivity of $\G$ that there is either a path from $r(e)$ to $r(f)$ (i.e., $r(e) \trianglerighteq r(f)$), or a path from $r(f)$ to $r(e)$ (i.e., $r(f) \trianglerighteq r(e)$), which contradicts the hypothesis that $v$ is a forked vertex.
\end{prf}\medskip

\begin{cor}
Let $\G$ be a directed graph and $\{\G_{i}\}_{i\in I}$ be the collection of all weakly connected components of $\G$, where $I$ is an index set. Then the congruence lattice $\mathscr{C}(\In)$ of $\In$ is distributive if and only if, for all $i\in I$, the congruence lattice $\mathscr{C}(\Ini)$ of $\Ini$ is distributive.
\end{cor}

\begin{prf} It is clear that $\G$ has no forked vertices if and only if $\G_i$ has no forked vertices, for all $i\in I$.
\end{prf}\medskip

\begin{eg}
As shown in Figure \ref{Figure 2}, the directed graph $\G_1$ has a forked vertex $v_1$. It has been proved in \cite[Example 2]{Luo} that the congruence lattice $\mathscr{C}(\Inx)$ is not lower semimodular.

For the directed graph $\G_2$ in Figure \ref{Figure 2}, it is easy to see that the graph inverse semigroup $\Iny =\{0, v_2, u_2, w_2, e_2, f_2, g_2, e_{2}^{*}, f_{2}^{*}, g_{2}^{*}\}$. There are 6 congruence triples in $\mathcal{CT}(\G_2)$: $T_1 = (\emptyset, \emptyset, \emptyset)$, $T_2 = (\{u_2\}, \emptyset, \emptyset)$, $T_3 = (\{w_2\}, \emptyset, \emptyset)$, $T_4 = (\{w_2\}, \{v_2\}, \emptyset\})$, $T_5 =(\{u_2, w_2\},\emptyset, \emptyset)$ and $T_6 = (\{u_2, v_2, w_2\}, \emptyset, \emptyset)$. It follows from the Hasse diagram of $\mathcal{CT}(\G_2)$ shown in Figure \ref{Figure 3} that $\mathscr{C}(\Iny)$ is distributive.
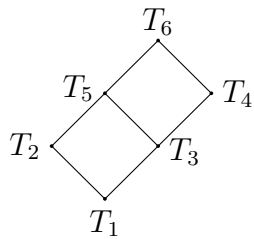
\begin{figure}[htbp]
\centering
\begin{tikzpicture}	[scale=0.7]
\fill (1,0) circle(1pt);
\fill (2,1) circle(1pt);
\fill (0,1) circle(1pt);
\fill (1,2) circle(1pt);
\fill (3,2) circle(1pt);
\fill (2,3) circle(1pt);

\node[below] at (1,0) {$T_1$};
\node[right] at (2,0.9) {$T_3$};
\node[right] at (3,2) {$T_4$};
\node[above] at (2,2.9) {$T_6$};
\node[left] at (1,2.1) {$T_5$};
\node[left] at (0,1) {$T_2$};

\draw  (1,0)--(2,1);
\draw  (2,1)--(3,2);
\draw  (3,2)--(2,3);
\draw  (1,2)--(2,3);
\draw  (1,2)--(0,1);
\draw  (1,2)--(2,1);
\draw  (1,0)--(0,1);
\end{tikzpicture}	
\caption{\footnotesize The Hasse diagram of the lattice $\mathcal{CT}(\G_2)$. \label{Figure 3}}
\end{figure}
\end{eg}

%%%%%%%%%%%%%%%%%%%%%%%%%%%%%%%%%%%%References%%%%%%%%%%%%%%%%%%%%%%%%%%%%%


\begin{thebibliography}{99}
\bibitem{AbramsPino} Abrams, G., Pino, G. A. (2005). The Leavitt path algebra of a graph. {\it J. Algebra} 293(2):319-334. \DOI{10.1016/j.jalgebra.2005.07.028}.

\bibitem{Anag} Anagnostopoulou-Merkouri, M., Mesyan, Z.,  Mitchell, J. D. (2022). Properties of congruence lattices of graph inverse semigroups. \arxiv{10.48550/arXiv.2108.08277}.

\bibitem{AshHall} Ash, C. J., Hall, T. E. (1975). Inverse semigroups on graphs. {\it Semigroup Forum}  11(1):140-145. \DOI{10.1007/BF02195262}.

\bibitem{Gratzer} Gr\"{a}tzer, G. (2011). Lattice Theory: Foundation. Springer, Birkh\"{a}user, Basel. \DOI{10.1007/978-3-0348-0018-1}.

\bibitem{Jones} Jones, P. R. (1983). On congruence lattices of regular semigroups. {\it J. Algebra} 82(1):18-39. \DOI{10.1016/0021-8693(83)90171-0}.

\bibitem{KumjianPaskRaeburn} Kumjian, A., Pask, D., Raeburn, I. (1998). Cuntz-Krieger algebras of directed graphs. {\it Pac. J. Math.} 184(1):161-174. \DOI{10.2140/pjm.1998.184.161}.

\bibitem{Luo} Luo, Y., Wang, Z. (2021). Semimodularity in congruence lattice of graph inverse semigroups. {\it Comm. Algebra} 49(6):2623--2632. \DOI{10.1080/00927872.2021.1879826}.

\bibitem{Meakin} Meakin, J., Milan, D., Wang, Z. (2021). On a class of inverse semigroups related to Leavitt path algebras. {\it Adv. in  Math} 384, 107729. \DOI{10.1016/j.aim.2021.107729}.

\bibitem{MesyanMitchell} Mesyan, Z., Mitchell, J. D. (2016). The structure of a graph inverse semigroup. {\it Semigroup Forum} 93(1):111-130. \DOI{10.1007/s00233-016-9793-x}.

\bibitem{Paterson} Paterson, A. L. T. (2002). Graph inverse semigroups, groupoids and their C$^\ast$-Algebras. {\it J. Operator theory} 48(3):645--662.

\bibitem{Wang} Wang, Z. (2019). Congruences on graph inverse semigroups. {\it J. Algebra} 534:51-64. \DOI{10.1016/j.jalgebra.2019.06.020}.
\end{thebibliography}
\end{document}